\documentclass[11pt]{amsart}   
\usepackage{amssymb}
\usepackage{amsthm}
\usepackage{amsmath}
\usepackage{amscd}

\newtheorem{thm}{Theorem}[section]
\newtheorem{prop}[thm]{Proposition}
\newtheorem{clly}[thm]{Corollary}

\theoremstyle{definition}
\newtheorem{defi}[thm]{Definition}

\newtheorem{nota}[thm]{Notation}

\theoremstyle{remark}
\newtheorem{rk}[thm]{Remark}

\newcommand{\ZZ}{{\mathbb{Z}}}
\newcommand{\RR}{{\mathbb{R}}}

\newcommand{\QQ}{{\mathbb{Q}}}

\newcommand{\TT}{{\mathbb{T}}}

\newcommand{\cstar}{\mbox{$C^*$}}

\newcommand{\gld}{\mbox{GL$_2({\mathbb Z})$}}
\newcommand{\glt}{\mbox{GL$_3({\mathbb Z})$}}

\newcommand{\gmn}{\mbox{$G_{\mu\nu}$}}
\newcommand{\gmnp}{\mbox{$G_{\mu'\nu'}$}}

\newcommand{\xau}{X^{\alpha,u}}
\newcommand{\xbus}{X^{\beta,u^*}} 

\newcommand{\supp}{\mbox{supp}}
\newcommand{\aut}{\mbox{Aut}}

\newcommand{\hcb}{\mbox{Hilbert \cstar-bimodule }}
\newcommand{\dcp}{D_{\mu'\nu'}^{c}}
\newcommand{\dc}{D_{\mu\nu}^{c}}

\title{Morita equivalence for Quantum Heisenberg Manifolds}
\author{Beatriz Abadie}
\address{Centro de Matem\'aticas. Facultad de Ciencias. 
Igu\'a 4225, CP 11 400, Montevideo, Uruguay.}
\email{abadie@cmat.edu.uy}
\thanks{Partially supported by Dinacyt (Proyecto Clemente Estable 8013), Uruguay.}

\date{\today}
\begin{document}

\begin{flushright}
\end{flushright}
\vspace{1cm}

\begin{abstract}

 We discuss Morita equivalence within the family $\{\dc: c\in \ZZ,\  c>0,\  \mu,\nu\in\RR\}$ of quantum Heisenberg manifolds. Morita equivalence classes are described in terms of the parameters $\mu$, $\nu$ and  the rank of the free abelian group $G_{\mu\nu}=2\mu\ZZ+2\nu\ZZ+\ZZ$ associated to the \cstar-algebra $\dc$.
 
\end{abstract}

\maketitle

{\bf Introduction.} Quantum Heisenberg manifolds $\{\dc: c\in \ZZ, c>0,\  \mu,\nu\in \RR\}$ were constructed by Rieffel in \cite{rfhm} as a quantization deformation of certain homogeneous spaces $H/N_c$, $H$ being the Heisenberg group.

It was shown in \cite[3.4]{fpa} that $K_0(\dc)\cong \ZZ^3\oplus\ZZ_c$,     which implies that $\dc$ and $D^{c'}_{\mu'\nu'}$ are not isomorphic unless $c=c'$.
Besides, $\dc$ and $D^c_{\mu'\nu'}$ are isomorphic when  $(2\mu,2\nu)$ and $(2\mu',2\nu')$ belong to the same orbit under the usual action of $\gld$ on $\TT^2$ ( \cite[Theorem 2.2]{pic}; see also \cite[3.3]{trace}). 
The range of traces on $\dc$ was discussed in \cite{trace}, where it was shown that the range of the homomorphism induced on $K_0(\dc)$ by any tracial state on $\dc$ has range $G_{\mu\nu}=\ZZ+2\mu \ZZ +2\nu\ZZ$. As a consequence (\cite[3.17]{trace}), the isomorphism condition  stated above turns out to be necessary when the rank of $G_{\mu\nu}$ is either 1 or 3. 
Rieffel showed in \cite{rfhm} that $\dc$ is simple if and only if $\{1,\mu,\nu\}$ is linearly independent over the field of rational numbers (i.e rank $G_{\mu\nu}=3$); it might be interesting to know whether in this case the classification can be made by means of the results of  Elliott and  Gong (\cite{eg}).

The quantum Heisenberg manifold $\dc$ was described in \cite{aee} as a crossed product by a Hilbert \cstar-bimodule. In order to discuss Morita equivalence within this family, we adapt to this setting some of the techniques employed in the analogous discussion for non-commutative tori (\cite{rfpm}) and Heisenberg \cstar-algebras (\cite{psm}). Thus  we generalize in Section 1  Green's result (discussed by Rieffel in \cite[Situation 10]{rfsit}) on the Morita equivalence of the crossed products $C_0(M/K)\rtimes H$ and $C_0(M/H)\rtimes K$, for free and proper commuting actions on a locally compact space $M$. This result provides the main tool  to discuss, in Section 2,   Morita equivalence for quantum Heisenberg manifolds.

\section{Morita equivalence of crossed products by certain Hilbert \cstar-bimodules over commutative \cstar-algebras.}

For a  Hilbert \cstar-bimodule $X$ over a \cstar-algebra $A$, the crossed product $A\rtimes X$ was introduced in \cite{aee} (see also \cite{pim}) as the universal \cstar-algebra for which there exist a *-homomorphism $i_A:A\longrightarrow A\rtimes X$
and a continuous linear map $i_X:X\longrightarrow A\rtimes X$ such that
\[ \begin{array}{ll}
i_X(ax)=i_A(a)i_X(x),\  i_A(\langle x,y\rangle_L)=i_X(x)i_X(y)^*.\\
i_X(xa)=i_X(x)i_A(a),\ i_A(\langle x,y\rangle_R)=i_X(x)^*i_X(y).
\end{array}\]
The crossed product $A\rtimes X$ carries a dual action $\delta$ of $S^1$, defined by $\delta_z(i_A(a))=i_A(a)$, $\delta_z(i_X(x))=zi_X(x)$, for $a\in A$, $x\in X$ and $z\in S^1$. Moreover, if a \cstar-algebra $B$ carries an action $\delta$ of $S^1$ such that $B$ is generated as a \cstar-algebra by the fixed point subalgebra $B_0=\{b\in B: \delta_z(b)=b\ \forall z\in S^1\}$ and the first spectral subspace $B_1=\{b\in B: \delta_z(b)=zb\ \forall z\in S^1\}$, then $B$ is isomorphic to $B_0\rtimes B_1$ (where $B_1$ has the obvious \hcb structure over $B_0$),  via an  isomorphism that takes the action $\delta$ into the dual action.

If $X$ is an $A$-\hcb  and $\alpha\in \aut(A)$, we denote by $X_\alpha$ the \hcb over $A$ obtained by leaving unchanged the left structure, and by setting
\[x\cdot_{X_{\alpha}} a:=x\alpha(a), \langle x, y\rangle^{X_\alpha}_R:=\alpha^{-1}(\langle x, y\rangle_R),\] where the undecorated notation refers to the original right structure of $X$. 

For $\alpha\in \aut(A)$ and the usual $A$-\hcb structure on $A$, the crossed product  $A\rtimes A_\alpha$ is easily checked to be the usual crossed product $A\rtimes_\alpha \ZZ$. 

\begin{defi}  Given a proper  action $\alpha$ of $\ZZ$ on a locally compact Hausdorff space $M$ and a unitary  $u\in C_b(M) $, let  $X^{\alpha,u}$ denote the set of functions $f\in C_b(M)$ satisfying $ f=u\alpha(f)$, and such that  the map  $x\mapsto |f(x)|$, which is constant on $\alpha$-orbits,  belongs to  $C_0(M/\alpha)$. 
Then $\xau$ is a Hilbert \cstar-bimodule over $C_0(M/\alpha)$ for pointwise multiplication on the left and the right, and  inner products given by $\langle f, g\rangle _L=f{\overline {g}}$, $\langle f,g\rangle _R={\overline {f}}g$.
\end{defi}

\begin{prop}
\label{unitary}
Let $\alpha$ and $\beta$  be free and proper commuting actions of $\ZZ$ on a locally compact Hausdorff space $M$, and let  $u$ be a unitary in $C_b(M)$. Then 
 the \cstar-algebras $C_0(M/\alpha)\rtimes  \xau_\beta$ and $C_0(M/\beta)\rtimes  \xbus_\alpha$ are Morita equivalent.
\end{prop}
\begin{proof}
Let $U:\ZZ\times \ZZ\longrightarrow {\mathcal U}(C_b(M))$ be given by 
\[U(n,k)=\left\{\begin{array}{ll}
1 & \mbox{ if either $n=0$ or $k=0$.}\\
\prod_{i\in S_k,j\in S_n}\alpha^i\beta^j(u), &\mbox{ for $n,k >0$}\\
\prod_{i\in S_k,j\in S_n}\alpha^i\beta^j(u^*), & \mbox{ for either $n$ or $k <0$, and $nk\neq 0$},\end{array}\right.\]
where $S_l=\{0,1,\dots l-1\}$ if $l> 0$ and $S_l=\{-1,-2,\dots l\}$ if $l<0$.
Straightforward computations show that $U(m+n,k)=U(m,k)\beta^m(U(n,k))$, and $U(n, k+l)=U(n,k)\alpha^k(U(n,l))$. 

Consider the proper actions    $\gamma^\alpha$ and  $\gamma^\beta$ of $\ZZ$ on $C_0(M)\rtimes_\beta \ZZ$ and  $C_0(M)\rtimes_\alpha \ZZ$, respectively, given by:
\[[\gamma^\alpha_k (\phi)](n)=U(n,k)\alpha^k[\phi(n)]\mbox { and } [\gamma^\beta_n (\psi)](k)=U^*(n,k) \beta^n[\psi(k)],\]
for $\phi\in C_c(\ZZ,C_0(M))\subset C_0(M)\rtimes_\beta \ZZ$ and $\psi\in C_c(\ZZ,C_0(M))\subset C_0(M)\rtimes_\alpha \ZZ$.

These two actions correspond, respectively, to $\gamma^{\alpha,U}$ and $\gamma^{\beta,U^*}$ in \cite[Propositions 1.2 and 2.1]{fpa}.
By virtue of \cite[Theorem 2.12]{fpa}, the generalized fixed-point algebras, in the sense of \cite[Definition 1.4]{rfpa},  $D^\alpha$ and $D^{\beta}$ of $C_0(M)\rtimes_\beta \ZZ$ and  $C_0(M)\rtimes_\alpha \ZZ$ under the actions $\gamma^\alpha$ and  $\gamma^\beta$, respectively, are Morita equivalent. The result will then be proved once we show that
$D^\alpha\cong C_0(M/\alpha)\rtimes  \xau_\beta$ and  $D^{\beta}\cong C_0(M/\beta)\rtimes  \xbus_\alpha$.

Recall from \cite[Proposition 2.1]{fpa} that $D^\alpha$ is defined to be the closed span in ${\mathcal M}(C_0(M)\rtimes_{\beta}\ZZ)$ of the set $\{P_\alpha(\phi^**\psi): \phi,\psi\in C_c(\ZZ\times M)\}$, where
\[P_{\alpha}(\phi)(x,n)=\sum_{k\in \ZZ}[\gamma^\alpha_k(\phi)](x,n),\]
for $\phi\in C_c(\ZZ\times M)\subset  C_0(M)\rtimes_{\beta}\ZZ,\  x\in M $, and $n\in \ZZ$. 

The \cstar-algebra $D^\alpha$ can also be described (\cite[Proposition 2.8]{fpa}) as the closure in  ${\mathcal M}(C_0(M)\rtimes_{\beta}\ZZ)$ of the $*$-subalgebra $C^{\alpha}=\{F\in C_c(\ZZ,C_b(M)): \gamma^\alpha(F)=F$   and $\pi_\alpha (\supp\ F(n))$ is precompact for all $n\in \ZZ\}$, 
where $\pi_{\alpha}$ denotes the canonical projection $\pi_{\alpha}:M\longrightarrow M/\alpha$.

Now, since $C^{\alpha}$ is contained in $ C_b(M)\rtimes_\beta \ZZ$, which is closed in $ {\mathcal M}(C_0(M)\rtimes_{\beta}\ZZ)$, so is $D^{\alpha}$. 
 Moreover, the \cstar-algebra $D^\alpha$ is invariant under the dual action $\hat{\beta}$ of $\TT$ on $C_b(M)\rtimes_\beta \ZZ$: 
\[\begin{array}{ll}
[\gamma^\alpha(\hat{\beta}_zF)](n,x)&=U(n,1)(x)(\hat{\beta}_z(F))(n,\alpha^{-1}x)\\
 &=U(n,1)(x)z^nF(n,\alpha^{-1}x)\\
&=z^nF(n,x)\\
&=(\hat{\beta}_zF)(n,x),
\end{array}\]
for $F\in C^\alpha$, $x\in M$ , $n\in \ZZ$,  and $z\in \TT.$ 
 Besides, $\supp\ (\hat{\beta}_z(F)(n))=\supp\ F(n)$ for all $n\in \ZZ$, so $\hat{\beta}_z(F)\in C^{\alpha}$ for all $z\in \TT$.

We next show that the action $\hat{\beta}$ on $D^\alpha$ is semi-saturated. That is, that, as a \cstar-algebra, $D^\alpha$ is generated by the fixed-point subalgebra  $D_0$ and the first spectral subspace $D_1=\{d\in D^\alpha: \hat{\beta}_z(d)=zd \ \forall z\in \TT\}$ for the restriction of the dual action $\hat{\beta}$. 

Since the maps $P_i:D^\alpha\longrightarrow D_i\ $ given by $P_i(a)=\int_{\TT}z^{-i}\hat{\beta}_z(a)\ dz$ are surjective contractions, $D_i$ is the closure of $P_i(C^\alpha)$.
Now, for $i=0,1$, $C_i=C^\alpha\cap F_i$, and $D_i=D^\alpha\cap F_i$, where $F_0$ and $F_1$ are, respectively, the fixed-point subalgebra and the first spectral subspace of $C_b(M)\rtimes_\beta\ZZ$, which are known to be the $\delta_i$-maps; that is, $F_i=\{F\in C_c(\ZZ,C_b(M)): \supp\ F=\{i\}\}$. 

Notice that
\[C^\alpha\cap F_0=\{f\delta_0:f\in C_b(M): \pi_\alpha(\supp\ f)\mbox{ is precompact and } f=\alpha(f)\}\]
 can be identified with $C_c(M/\alpha)$ via $f\delta_0\mapsto \tilde{f}$, where $\tilde{f}\circ \pi_\alpha=f$, and that this map extends to a $*$-isomorphism between $D_0$ and $C_0(M/\alpha)$.

Now, $D_1$ is a \hcb over $D_0$ for  
\begin{subequations}
\label{d}
\begin{gather}
(f\delta_0)*(g\delta_1)=(fg)\delta_1\mbox{ and } (g\delta_1)*(f\delta_0)=(g\beta(f))\delta_1,\\
 \langle f\delta_1,g\delta_1\rangle_L = (f\delta_1)*(g\delta_1)^*=(f\bar{g})\delta_0,\\ 
 \langle f\delta_1,g\delta_1\rangle_R=(f\delta_1)^**(g\delta_1)=(\beta^{-1}(\bar{f}g))\delta_0 .
\end{gather}
\end{subequations}
Notice  that $D_1$ is full on the left (and on the right, by a similar argument) as a \hcb over $C_0(M/\alpha)$. For $<D_1,D_1>_L$, the closed linear span in $C_0(M/\alpha)\cong D_0$ of the set $\{\langle\ f\delta_1, g\delta_1\rangle_L: f\delta_1,g\delta_1\in D_1\}$, is a closed ideal of $C_0(M/\alpha)$. 
Therefore, unless $<D_1,D_1>_L=C_0(M/\alpha)$, there exists $x_0\in M$ such that $f(x_0)=0$ for all $f\delta_1\in D_1$. 

Now, given $x_0\in M$,  we can choose (\cite[Situation 10]{rfsit}) a neighborhood $U$ of $x_0$ such that $U\cap \alpha^k(U)=\emptyset$ for  $k\neq 0$. 
Let $g\in C_c(M)^+$ be such that $\supp\ g\subset U$ and $g(x_0)=1$. 

Then 
\begin{multline}\notag
[P_\alpha\big((g^{1/2}\delta_0)^**(g^{1/2}\delta_1)\big)](x,n)=\\=(P_\alpha(g\delta_1))(x,n)=\big(\sum_k U(1,k)(x)g(\alpha^{-k}(x))\big)\delta_1(n),
\end{multline}
so  $P_\alpha\big((g^{1/2}\delta_0)^**(g^{1/2}\delta_1)\big)\in D_1$, and equals 1 at $(x_0,1)$.

In order to prove that $C^\alpha\subset C^*(D_0,D_1)$, it suffices to show that $f\delta_k\in C^*(D_0,D_1)$ for $f\delta_k\in C^\alpha$, $k\in \ZZ$ ; since $C^\alpha$ is closed under involution, we may assume that $k\geq 0$. 
We show this fact, which clearly holds for $k=0$ and $k=1$, by induction on $k$. 

If $f\delta_k\in C^\alpha$ and $\epsilon >0$, since $\pi_\alpha(\supp \ f)$ is precompact in $M/\alpha$, and $D_1$ is full over $C_0(M/\alpha)$, we can find $\phi_i,\psi_i\in D_1, i=1,\dots ,p,$ such that $\|\sum_i (\phi_i*\psi^*_i)*f\delta_1-f\delta_1\|_{D^\alpha}=\|\sum_i\langle \phi_i,\psi_i\rangle_L f-f\|_{C_b(M)}<\epsilon $. 

Now, since $\phi_i$ and  $\psi_i^**f$ belong to $C^*(D_0,D_1)$ for $i=1,\dots,p$, so does $f$.
This shows that $D^\alpha =C^*(D_0,D_1)$ and, consequently,  by \cite[Theorem 3.1]{aee}, that $D^\alpha\cong D_0\rtimes D_1$.

It only remains to notice now that $D_0\rtimes D_1\cong C_0(M/\alpha)\rtimes X^{\alpha,u}_\beta$.
As noticed above, $D_0$ is isomorphic to  $C_0(M/\alpha)$. 
On the other hand, the map $f\delta_1\mapsto f$ takes $C^\alpha\cap F_1$ to $X^{\alpha,u}_\beta$. By keeping track of the formulae in \ref{d}, one easily checks that that map is an isometry, so it extends to an isometry from $D_1$ to $X^{\alpha,u}_\beta$, which is onto because its image  contains the dense set:
\[X_0^{\alpha,u}=\{f\in X^{\alpha,u}_\beta \mbox{:  the map } x\mapsto |f(x)| \mbox{ is compactly supported on $M/\alpha$}\}.\]
(Notice that $X_0^{\alpha,u}$ is dense in $X^{\alpha,u}_\beta$, because, if $\{e_\lambda\}$ is an approximate identity for $C_c(M/\alpha)$, then $e_\lambda f$ converges to $f$ for all $f\in X^{\alpha,u}_\beta$).

This shows that $D^\alpha$ is isomorphic to $C_0(M/\alpha)\rtimes X^{\alpha,u}_\beta$. Analogously, $D^\beta$ is isomorphic to $C_0(M/\beta)\rtimes X^{\beta,u^*}_\alpha$.

\end{proof}

\section{Morita equivalence for quantum Heisenberg manifolds}

In \cite{aee} (see also \cite[2]{pic}) the quantum Heisenberg manifold $\dc$ was shown to be the crossed product of $C(\TT^2)$, the \cstar-algebra of continuous functions on the torus, by  the \hcb 
$M^c_{\alpha_{\mu\nu}}$,   where $\alpha_{\mu\nu}(x,y)=(x+2\mu, y+2\nu)$,  and \[M^c=\{f\in C_b(\RR\times\TT): f(x+1,y)=e^{-2\pi icy}f(x,y)\}\] 
is the Hilbert \cstar-bimodule obtained by letting  $C(\TT^2)$ act  by pointwise product, and by defining the inner products $\langle f,g\rangle _L=f\overline{g}$, $\langle f,g\rangle _R=\overline{f}g$.

\begin{rk}
\label{sorb}
The \cstar-algebras $\dc$ and $\dcp$  are isomorphic when the projections of $(2\mu,2\nu)$ and $(2\mu',2\nu')$ on the torus are in the same orbit under the usual action of $\gld$ (\cite[Theorem 2.2]{pic}, see also \cite[Remark 3.3]{trace}).
\end{rk}
\begin{prop}
\label{fmu}
Let $\mu\neq 0$. Then $\dc$ and $D^c_{\frac{1}{4\mu},\frac{\nu}{2\mu}}$ are Morita equivalent.
\end{prop}

\begin{proof}
We follow the lines of \cite[1.1]{rfpm}  and apply  Proposition \ref{unitary} to the following setting:  $\alpha$ and $\beta$ consist of  translation on $\RR\times \TT$ by $(\frac{1}{2\mu},0)$ and $(1,2\nu)$, respectively, and  $u\in C_b(\RR\times \TT)$ is given by $u(x,y)=e(-cy)$, where $\TT$ is viewed as $\RR/\ZZ$ and, for a real number $h$, $e(h)=e^{2\pi ih}$. 

Then, by Proposition \ref{unitary},   $C((\RR\times\TT)/\alpha)\rtimes X^{\alpha,u}_\beta$ and 
 $C((\RR\times\TT)/\beta)\rtimes X^{\beta,u^*}_\alpha$ are Morita equivalent, where 
\begin{multline}
\notag
X^{\alpha,u}=\{F\in C_b(\RR\times \TT):\ F(x-\frac{1}{2\mu},y)=e(cy)F(x,y)\}\mbox { and }\\
X^{\beta,u^*}=\{F\in C_b(\RR\times \TT):\ F(x-1,y-2\nu)=e(-cy)F(x,y)\}\\
=\{F\in C_b(\RR\times \TT):\ F(x+1,y+2\nu)=e(c(y+2\nu))F(x,y)\} .
\end{multline}
 Let $H_{\alpha}:C(\TT^2)\longrightarrow C((\RR\times\TT)/\alpha)$ and  $H_{\beta}:C(\TT^2)\longrightarrow C((\RR\times\TT)/\beta)$ be the isomorphisms given by:
\[(H_\alpha\phi)(x,y)=\phi(2\mu x,y),\ (H_\beta\phi)(x,y)=\phi(x,2\nu x-y),\]
and, for $(\mu',\nu')=(\frac{1}{4\mu},\frac{\nu}{2\mu})$,  set 
\[J_\alpha: M^c_{\alpha_{\mu\nu}}\longrightarrow X^{\alpha,u}_\beta\mbox{  and  }J_\beta:  M^c_{\alpha_{\mu'\nu'}}\longrightarrow X^{\beta,u^*}_\alpha,\]
\[(J_{\alpha}f)(x,y)=f(2\mu x,y),\ (J_{\beta}f)(x,y)=e(cx(x+1)\nu)f(x,2\nu x-y).\]

Notice that 
\[(J_{\alpha}f)(x-\frac{1}{2\mu},y)=f(2\mu x-1,y) =e(cy)(J_{\alpha}f)(x,y),\]
and 
\[\begin{array}{ll}
(J_{\beta}f)(x+1,y+2\nu)&=e(c(x+1)(x+2)\nu)f(x+1,2\nu x-y)\\
&= e(c(x+1)(x+2)\nu)e(-c(2\nu x-y))f(x,2\nu x-y)\\ 
&=e(c(y+2\nu))(J_{\beta}f)(x,y),
\end{array}\]
so the definitions make sense. 

For $i=\alpha,\beta$, it is easily checked that $J_i$ is a bijection and that, for $\phi\in C(\TT^2)$, $f,g\in M^c$:
\[ J_i(\phi\cdot f)=H_i(\phi)\cdot J_i(f),\ J_i(f \cdot \phi)=J_i(f)\cdot H_i(\phi),\]
\[ \langle J_if,J_ig\rangle_L=H_i(\langle f,g\rangle_L),\  \langle J_if,J_ig\rangle_R=H_i(\langle f,g\rangle_R).\]
This shows that $\dc=C(\TT^2)\rtimes M^c_{\alpha_{\mu\nu}}$ and $\dcp=C(\TT^2)\rtimes M^c_{\alpha_{\mu'\nu'}}$ are isomorphic, respectively, to $C((\RR\times\TT)/\alpha)\rtimes X^{\alpha,u}_\beta$ and 
 $C((\RR\times\TT)/\beta)\rtimes X^{\beta,u^*}_\alpha$, and they are, consequently, Morita equivalent to each other.
\end{proof}

\begin{clly}
\label{gl2}
Let $\mu\not\in \QQ$,  and  let $A= \left(\begin{array}{ll}
a & b \\
c & d 
\end{array}
\right) \in \gld$. If 
\[ 2\mu'=\frac {2a\mu+b}{2c\mu+d} \mbox{ and } 2\nu'=\frac {2\nu}{2c\mu+d}, \]
then the quantum Heisenberg  manifolds $\dc$ and $\dcp$ are Morita equivalent.
\end{clly}
\begin{proof}

It suffices to check the statement for  $ A_1= \left(\begin{array}{ll}
1 &  1\\
0 &  1
\end{array}
\right)$ and $A_2 =\left(\begin{array}{ll}
0 &  1\\
1 &  0
\end{array}
\right)$
, since $A_1$ and $A_2$ generate $\gld$ (\cite[Appendix B]{kur}), and $(\mu,\nu)\mapsto (\mu',\nu')$ defines an action of $\gld$ on $\big(\RR\setminus\QQ\big)\times \RR$. 
 For $A=A_1$ we get isomorphic  \cstar-algebras by Remark \ref{sorb}. For $A=A_2$, we get $(\mu',\nu')=(\frac{1}{4\mu},\frac{\nu}{2\mu})$, and the result follows from Proposition \ref{fmu}.

\end{proof}

\begin{prop}
\label{gl3}
Let $\{1,\mu,\nu\}$ be linearly independent over $\QQ$, and  let $A= \left(\begin{array}{lll}
a & b & c\\
d & e & f \\
g & h & i
\end{array}
\right) \in \glt$. If 
\[ 2\mu'=\frac {2a\mu+2b\nu +c}{2g\mu+2h\nu+i} \mbox{ and } 2\nu'= \frac {2d\mu+2e\nu +f}{2g\mu+2h\nu+i}, \]
then the quantum Heisenberg  manifolds $\dc$ and $\dcp$ are Morita equivalent.
\end{prop}

\begin{proof}
As in the proof of Theorem 1.7 in \cite{psm}, $ A=A_1A_2A_3$, where
\[A_1=\left(\begin{array}{lll}
A & B & C\\
D & E & F \\
0 & 0 & 1
\end{array}
\right),\  A_2  \left(\begin{array}{lll}
G & 0 & H\\
0 & 1 & 0 \\
I & 0 & J
\end{array}
\right),\ A_3= \left(\begin{array}{lll}
K & L & 0\\
M & N & 0 \\
0 & 0 & 1
\end{array}
\right),\]
and $A_i\in\glt$, for $i=1,2,3$.

Since the map $(\mu,\nu)\mapsto (\mu',\nu')$ defines an action of $\glt$ on the set $\{(\mu,\nu)\in\RR^2: \{1,\mu,\nu\} \mbox{ is linearly independent over $\QQ$}\}$, it suffices to check the statement for $A_i$, $i=1,2,3.$

For $A=A_1$ and $A=A_3$ the \cstar-algebras $\dc$ and $\dcp$ are isomorphic by Remark \ref{sorb}. Thus it suffices to show the result for $A=A_2$. 
The map $ \left(\begin{array}{ll}
G &  H\\
I &  J
\end{array}
\right) \mapsto  \left(\begin{array}{lll}
G & 0 & H\\
0 & 1 & 0 \\
I & 0 & J
\end{array}
\right)$ is a group homomorphism from  $\gld$ into $\glt$, and $\gld$ is generated by  
$\left(\begin{array}{ll}
1 &  1\\
0 &  1
\end{array}
\right)$  and 
$\left(\begin{array}{ll}
0 &  1\\
1 &  0
\end{array}
\right)$, which implies that we only need to prove the statement for 
$A_1= \left(\begin{array}{lll}
1 & 0 & 1\\
0 & 1 & 0 \\
0 & 0 & 1
\end{array}
\right)$ and $ A_2= \left(\begin{array}{lll}
0 & 0 & 1\\
0 & 1 & 0 \\
1 & 0 & 0
\end{array}
\right)$. For $A_1$ we get $2\mu'=2\mu+1$, $2\nu'=2\nu$, so  $\dc$ and $\dcp$ are isomorphic  by Remark \ref{sorb}. Proposition \ref{fmu} takes care of the case  $A=A_2$, since then  we have $(\mu',\nu')=(\frac{1}{4\mu},\frac{\nu}{2\mu})$.
\end{proof}
\begin{nota}
We denote by $\gmn$ the  subgroup of $\RR$ generated by $\{1,2\mu,2\nu\}$. It was shown in \cite[Theorem 3.16]{trace} that the homomorphism induced on $K_0(\dc)$ by any tracial state on $\dc$ has range $\gmn$.  
\end{nota}

\begin{rk}
\label{rank12}
 If rank $\gmn=2$, then there exist an irrational number $\nu'$ and  integers $p,q\in \ZZ,\  p\neq 0$, $(p,q)=1$,  such that  $\dc$ and $D^c_{\frac{p}{2q},\nu'}$ are isomorphic. 
\end{rk}
\begin{proof}
We proceed as in \cite[Proposition 1.5]{pcs}. Let $\mu_0=2\mu$, $\nu_0=2\nu$. Since the group generated by $\{1,\mu_0,\nu_0\}$ has rank 2, either $\mu_0$ or $\nu_0$ is irrational. We may assume that $\nu_0$ is irrational, because, by Remark \ref{sorb}, $\dc$ and $D^c_{\nu\mu}$ are isomorphic. 
 Besides, there exist $M,N,P\in \ZZ$, with $N\neq 0$ such that $M+N\mu_0+P\nu_0=0$, so we have $\mu_0=\frac{k}{l}\nu_0+\frac{m}{n}$, with $(k,l)=1$. If $k=0$, then $\mu_0\in \QQ$, and we are done. Otherwise take $a,b\in \ZZ$ such that $ak+bl=1$, so that $ \left(\begin{array}{ll}
-l & k\\
a & b
\end{array}\right) \in\gld$, and set 
\[(\mu'_0,\nu'_0)=\left(\begin{array}{ll}
-l & k\\
a & b
\end{array}\right)(\mu_0,\nu_0).\]
 Then 
\[\mu'_0=-l(\frac{k}{l}\nu_0+\frac{m}{n})+k\nu_0=\frac{-lm}{n}\in \QQ.\]
and
\[\nu'_0=a(\frac{k}{l}\nu_0+\frac{m}{n})+b\nu_0=\frac{1}{l}\nu_0+\frac{am}{n}\not\in\QQ,\]
We now take $\nu'=\nu'_0/2$ and $p/q=\mu'_0$,  in lowest terms. By Remark \ref{sorb} $\dc$ and $D^c_{\frac{p}{2q},\nu'}$ are isomorphic.
\end{proof}

\begin{prop}
\label{dif}
Let $p$ and $q$ be non-zero integers such that $(p,q)=1$, and let $\nu\in\RR$. Then $D^c_{\frac{p}{2q},\nu}$ is Morita equivalent to $D^c_{0,q\nu}$.
\end{prop}
\begin{proof}

By Remark \ref{sorb} we may assume that $p$ and $q$ are positive. By applying Proposition \ref{fmu} to $(\mu,\nu)=(q/2,\nu)$,  we get that $D^c_{0,\nu}\cong D^c_{q/2,\nu}$ is Morita equivalent to $D^c_{\frac{1}{2q},\frac{\nu}{q}}$, thus proving the proposition for $p=1$.
For $p>1$, let $r_0=q$, $r_1=p$ , and, if $r_{i+1}\neq 1$, define $r_{i+2}$ by $r_i=m_{i+1}r_{i+1} +r_{i+2}$, where $0\leq r_{i+2}<r_{i+1}$, and $m_{i+1}\in \ZZ$.

  Actually, $r_{i+2}>0$: otherwise  $r_{i+1}$ divides $r_i$, and it follows that $r_{i+1}$  divides $r_{j}$ for all $j\leq i$. In particular, $r_{i+1}$ divides both $p$ and $q$, which contradicts the fact that $r_{i+1}\neq 1$. 
Now, since $r_{i+1}<r_i$, there is an index $i_0$ for which $r_{i_0}=1$.

  On the other hand, it follows from Proposition \ref{fmu} that, for any real number $\kappa$, $D^c_{\frac{r_i}{2r_{i-1}},\kappa}$ is Morita equivalent to $D^c_{\frac{r_{i-1}}{2r_i},\kappa\frac{r_{i-1}}{r_i}}$, which in turn is isomorphic to $D^c_{\frac{r_{i+1}}{2r_i},\kappa\frac{ r_{i-1}}{r_i}}$. 

Thus $D^c_{\frac{p}{2q},\nu}=D^c_{\frac{r_1}{2r_0},\nu}$ is Morita equivalent to $D^c_{\frac{r_j}{2r_{j-1}}, \frac{q\nu}{r_{j-1}}}$ for any $j\leq i_0$. In particular, for $j=i_0$, we have that $D^c_{\frac{p}{2q},\nu}$ is Morita equivalent to $D^c_{\frac{1}{2r_{i_0-1}},\frac{\nu q}{r_{i_0-1}}}$, which, as shown above, is Morita equivalent to $D^c_{0,\nu q}$.
\end{proof}

\begin{thm}

Two quantum Heisenberg manifolds $\dc$, $D^{c'}_{\mu'\nu'}$ are Morita equivalent if and only if $c=c'$ and there exists a positive real number $r$ such that
\[\ZZ+2\mu\ZZ+2\nu\ZZ=r(\ZZ+2\mu'\ZZ+2\nu'\ZZ ).\]
In particular, the rank of the free abelian group $\gmn=\ZZ+2\mu\ZZ+2\nu\ZZ$ is the same for Morita equivalent quantum Heisenberg manifolds, and:
\begin{enumerate}
\item  If rank $\gmn =1= $ rank $\gmnp$, then $\dc$ is Morita equivalent to $D^{c}_{\mu',\nu'}$. In particular, $\dc$ is Morita equivalent to the commutative Heisenberg manifold $D^c_{0,0}$.

\item  If rank $\gmn =2=$ rank $\gmnp$, let $\{\alpha,\frac{p}{q}\}$ and $\{\alpha',\frac{p'}{q'}\}$ be bases of $\gmn$ and $\gmnp$, respectively, where $\alpha$ and $\alpha'$ are irrational numbers and  $p,p',q,q'\in\ZZ$, $(p,q)=(p'q')=1$. Then  $D^c_{\mu,\nu}$ and $D^c_{\mu',\nu'}$ are Morita equivalent if and only if there exists
$ \left(\begin{array}{ll}
a & b\\
c & d
\end{array}\right) \in\gld$ such that 
\[q'\alpha '=\frac{aq\alpha '+b}{cq\alpha '+d}.\]
In particular, $D^c_{\frac{p}{2q},\nu}$ is Morita equivalent to $D^c_{0,q\nu}$.
\item  If rank $\gmn =3=$ rank $\gmnp$, then $\dc$ and $\dcp$ are Morita equivalent if and only if there exists $\left(\begin{array}{lll}
a & b & c\\
d & e & f \\
g & h & i
\end{array}
\right) \in \glt$ such that
\[ 2\mu'= \frac {2a\mu+2b\nu +c}{2g\mu+2h\nu+i} \mbox{ and }2 \nu'= \frac {2d\mu+2e\nu +f}{2g\mu+2h\nu+i}.\]

\end{enumerate}

\end{thm}
\begin{proof}

It was shown in \cite[3.4]{fpa} that $K_0(\dc)=\ZZ^3\oplus \ZZ_c$, which implies that $\dc$ and $D^{c'}_{\mu'\nu'}$ are not Morita equivalent for $c\neq c'$.

  Besides (\cite[Theorem 3.16]{trace}), all tracial states on $\dc$ induce the same homomorphism on $K_0(\dc)$, whose range is the group $G_{\mu\nu}=2\mu \ZZ+2\nu\ZZ+\ZZ$. Since (\cite[2.2]{rfir}) there is a bijection between finite traces of Morita equivalent unital \cstar-algebras we must have  $G_{\mu\nu}=rG_{\mu'\nu'}$ for some positive real number $r$ when $\dc$ and $\dcp$ are Morita equivalent.  An immediate consequence of this fact is that the rank of $G_{\mu\nu}$ is  invariant under Morita equivalence.
 
If rank $\gmn=1$, then, by \cite[Remark 3.5]{trace},  $\dc$ is isomorphic to $D^c_{0,\frac{1}{2p}}$ for some non-zero integer $p$, so $\dc$  is isomorphic to $D^c_{\frac{1}{2p},0}$ by Remark \ref{sorb}. Now, by Proposition \ref{dif}, $D^c_{\frac{1}{2p},0}$ is Morita equivalent  $D^c_{0,0}$.

If rank $\gmn=2=$ rank $\gmnp$ and $G_{\mu\nu}=rG_{\mu'\nu'}$ for some positive $r$, let  $\{\alpha,\frac{p}{q}\}$ and $\{\alpha',\frac{p'}{q'}\}$ be bases of $\gmn$ and $\gmnp$, respectively, where $\alpha,\ \alpha'$ are irrational numbers, and $p,p',q,q'$ are integers, with $(p,q)=(p',q')=1$. Since $\ZZ\subset \gmn\ (\gmnp)$ we have that $p\  (p')=\pm1$ and, by Remark \ref{sorb}, we may assume $p=p'=1$. Then  we have that $\alpha \ZZ+1/q\ZZ=r(\alpha' \ZZ+1/q'\ZZ)$, which implies that $\alpha q\ZZ+\ZZ=(rq/q')(\alpha'q' \ZZ+\ZZ)$. A standard  argument shows that 
\[q\alpha=\frac {aq'\alpha'+b}{cq'\alpha'+d}\mbox{ for some } \left(\begin{array}{ll}
a & b\\
c & d
\end{array}\right) \in\gld.\]
Therefore $D^c_{\frac{q\alpha}{2},0}$ and $D^c_{\frac{q'\alpha'}{2},0}$ are Morita equivalent by Corollary \ref{gl2}. 

On the other hand, by Remark \ref{rank12}, $\dc$ and $D^c_{\mu',\nu'}$ are isomorphic, respectively, to $D^c_{\frac{m}{2n},\beta}$ and  $D^c_{\frac{m'}{2n'},\beta'}$, for some irrational numbers $\beta$ and $\beta'$ and integers $m,m',n,n'$ such that $(m,n)=(m',n')=1$.
Therefore  $\{2\beta,\frac{1}{n}\}$ and $\{2\beta',\frac{1}{n'}\}$ are bases of $\gmn$ and $G_{\mu'\nu'}$, respectively, and it follows from the argument above that $D^c_{n\beta,0}$ and $D^c_{n'\beta',0}$ are Morita equivalent. It only remains to notice now that, by Proposition \ref{dif} and Remark \ref{sorb}, $D^c_{n\beta,0}$ and $D^c_{n'\beta',0}$ are Morita equivalent to $\dc$ and $D^c_{\mu'\nu'}$, respectively.

Finally,  if rank $\gmn=3=$ rank $\gmnp$ and $\gmn=r\gmnp$ for some positive $r$, then let  $A=\left(\begin{array}{lll}
a & b & c\\
d & e & f \\
g & h & i
\end{array}
\right) \in \glt$ be the transpose of the matrix that changes coordinates between the bases  $\{2r\mu',2r\nu',r\}$ and  $\{2\mu,2\nu,1\}$ of $\gmn$. Then
\[ 2\mu'=\frac {2a\mu+2b\nu +c}{2g\mu+2h\nu+i} \mbox{ and }2 \nu'= \frac {2d\mu+2e\nu +f}{2g\mu+2h\nu+i},\]
which implies, by Proposition \ref{gl3}, that $\dc$ and $\dcp$ are Morita equivalent.

\end{proof}

\noindent

\end{document}